\documentclass[a4paper,10pt,reqno]{amsart}
\usepackage[dvipsnames]{xcolor}
\usepackage[T1]{fontenc}
\usepackage{hyperref}
\usepackage{bm}   

\hypersetup{
colorlinks = True,
linkcolor=Bittersweet,
citecolor=blue,
}
\definecolor{mycolor}{RGB}{226, 90, 104}
\definecolor{bluey}{rgb}{0.0, 0.0, 1.0}
\definecolor{greeny}{rgb}{0.0, 0.5, 0.0}
\definecolor{purpley}{rgb}{0.5, 0.0, 0.5}
\definecolor{orangey}{rgb}{1.0, 0.65, 0.0}
\definecolor{yellowy}{rgb}{1.0, 1.0, 0.0}
\definecolor{limegreeny}{rgb}{0.2, 0.8, 0.2}

\usepackage[nameinlink,capitalize, noabbrev]{cleveref}
\usepackage{multirow}
\usepackage{mathtools}
\usepackage{amsfonts,amsmath,mathrsfs,amssymb,amsthm,amscd,latexsym,amstext,amsxtra, mathabx}
\usepackage{etoolbox}
\usepackage[all]{xy}
\usepackage{enumerate}
\usepackage{tabularx}
\xyoption{all}
\usepackage{srcltx} 
\usepackage[english]{babel}
\usepackage{epsfig}
\usepackage[]{epigraph}
\usepackage{xy}
\usepackage{enumerate}
\usepackage{graphicx}
\usepackage{float}
\usepackage{cancel}
\usepackage[numbers]{natbib}

\DeclarePairedDelimiter{\floor}{\lfloor}{\rfloor}
\DeclarePairedDelimiter{\bracket}{\langle}{\rangle}

\newtheorem{theorem}{Theorem}[section]
\newtheorem{lemma}[theorem]{Lemma}
\newtheorem{corollary}[theorem]{Corollary}
\newtheorem{prop}[theorem]{Proposition}
\newtheorem{remark}[theorem]{Remark}

\newtheorem{definition}[theorem]{Definition}

\newtheorem{example}[theorem]{Example}

\newcommand{\nc}{\newcommand}
\nc{\cH}{{\mathcal H}}
\nc{\cA}{{\mathcal A}}
\nc{\cG}{{\mathcal G}}
\nc{\cC}{{\mathcal C}}
\nc{\cD}{{\mathcal D}}
\nc{\cO}{{\mathcal O}}
\nc{\cI}{{\mathcal I}}
\nc{\cB}{{\mathcal B}}
\nc{\cY}{{\mathcal Y}}
\nc{\cK}{{\mathcal K}}
\nc{\cX}{{\mathcal X}}
\nc{\cS}{{\mathcal S}}
\nc{\cE}{{\mathcal E}}
\nc{\cF}{{\mathcal F}}
\nc{\cZ}{{\mathcal Z}}
\nc{\cQ}{{\mathcal Q}}
\nc{\cN}{{\mathcal N}}
\nc{\cP}{{\mathcal P}}
\nc{\cL}{{\mathcal L}}
\nc{\cM}{{\mathcal M}}
\nc{\cT}{{\mathcal T}}
\nc{\cW}{{\mathcal W}}
\nc{\cU}{{\mathcal U}}
\nc{\cJ}{{\mathcal J}}
\nc{\cV}{{\mathcal V}}
\nc{\cR}{{\mathcal R}}
\nc{\bH}{{\mathbb H}}
\nc{\bA}{{\mathbb A}}
\nc{\bG}{{\mathbb G}}
\nc{\bC}{{\mathbb C}}
\nc{\bO}{{\mathbb O}}
\nc{\bI}{{\mathbb I}}
\nc{\bB}{{\mathbb B}}
\nc{\bY}{{\mathbb Y}}
\nc{\bK}{{\mathbb K}}
\nc{\bX}{{\mathbb X}}
\nc{\bS}{{\mathbb S}}
\nc{\bE}{{\mathbb E}}
\nc{\bF}{{\mathbb F}}
\nc{\bZ}{{\mathbb Z}}
\nc{\bQ}{{\mathbb Q}}
\nc{\bN}{{\mathbb N}}
\nc{\bP}{{\mathbb P}}
\nc{\bL}{{\mathbb L}}
\nc{\bM}{{\mathbb M}}
\nc{\bT}{{\mathbb T}}
\nc{\bW}{{\mathbb W}}
\nc{\bU}{{\mathbb U}}
\nc{\bD}{{\mathbb D}}
\nc{\bJ}{{\mathbb J}}
\nc{\bV}{{\mathbb V}}
\nc{\bbZ}{{\mathbb Z}}
\nc{\bR}{{\mathbb R}}
\nc{\fX}{{\mathfrak X}}
\nc{\Ida}{{\mathfrak a}}
\nc{\Idb}{{\mathfrak b}}
\nc{\Idp}{{\mathfrak p}}
\nc{\Idm}{{\mathfrak m}}
\nc{\fr}{{\rightarrow}}
\nc{\co}{{\nabla}}
\nc{\ok}{{\overline{K}}}
\nc{\Aut}{{\mbox{Aut}}}
\nc{\Int}{{\mbox{Int}}}
\nc{\la}{\longrightarrow}
\nc{\elemen}{\mathcal{E}(\bR^d)}
\nc{\Spec}{\mbox{Spec}}
\nc{\Proj}{\mbox{Proj}}
\nc{\Sym}{{\mathcal{S}}ym}
\nc{\divi}{\mbox{div}}
\nc{\Divi}{\mbox{Div}}
\nc{\Gal}{\mbox{Gal}}
\nc{\Bl}{\mbox{Bl}}
\nc{\GL}{\mbox{GL}}
\nc{\PGL}{\mbox{PGL}}
\nc{\mult}{\mbox{mult}}
\nc{\Conv}{\mbox{Conv}}
\nc{\Supp}{\mbox{Supp}}
\nc{\Pic}{\mbox{Pic}}
\nc{\cu}{{\overlineline{\nabla}}}
\nc{\con}{{\nabla}}
\nc{\vol}{{\mbox{vol}}}
\nc{\logc}{c}
\nc{\Ch}{\mbox{Ch}}
\nc{\Td}{\mbox{Td}}
\nc{\Cos}{\mbox{Cos}}
\nc{\Sin}{\mbox{Sin}}
\nc{\re}{\mbox{Re}}
\nc{\im}{\mbox{Im}}
\nc{\om}{\overline{m}}
\title[Reciprocity For Dedekind Sums via Conical Zeta Values]{Reciprocity For Dedekind Sums via Conical Zeta Values}
\author[Yerko Torres-Nova]{Yerko Torres-Nova}
\email{yerko.torres@usm.cl}
\address{Departamento de Matem\'atica, Universidad T\'ecnica Federico Santa Mar\'ia, Valpara\'iso, Chile.}

\begin{document}

\begin{abstract}
We study reciprocity formulas for Dedekind sums associated with absolutely continuous functions, extending the classical Dedekind-Rademacher reciprocity formula. In particular, we treat the case of periodic Bernoulli functions. Our approach generalizes an integral method and uses Fourier analysis to show that the reciprocity for polynomial-type functions admits a geometric interpretation in terms of conical zeta values.
\end{abstract}

\maketitle


\tableofcontents

\section{Introduction}

Given positive integers $\nu_0,\nu_1,\nu_2$, the classical Dedekind sum is defined as
$$S(\nu_1,\nu_2 | \nu_0) := \sum_{i=1}^{\nu_0-1}  b_1\left( \frac{i\nu_1 }{\nu_0} \right) b_1\left( \frac{i\nu_2}{\nu_0} \right),$$
where $b_1(x)$ is the \emph{first periodic Bernoulli function}, i.e.,
$$b_1:\bR \to \bR,\quad b_1(x) = 
    \{x\}-1/2,$$
and $\{x \} = x-\floor{x}$ denotes the fractional part of $x$. The Dedekind-Rademacher reciprocity theorem \cite{rademacherreci,rademacheremil} asserts that if $\nu_0,\nu_1,\nu_2$ are pairwise coprime positive integers, then
\begin{equation}\label{rademacherreci}
   \cR_{b_1}(\bm{\nu}) := S(\nu_1,\nu_2|\nu_0) + S(\nu_0,\nu_2|\nu_1) + S(\nu_1,\nu_0|\nu_2) = \frac{\nu_0^2 + \nu_1^2+\nu_2^2}{12\nu_0 \nu_1\nu_2} - \frac{1}{4}.
\end{equation}
This reciprocity law is fundamental in many areas and provides techniques that link arithmetic, geometry, and topology (see, for instance, \cite{2009dedekindCfbirs}). Due to its significance, numerous approaches have been developed to study its nature and to derive generalizations in various directions \cite{carlitz74,zagier1973}. For our purposes, a key example arises by replacing the function $b_1(x)$ with higher-degree periodic Bernoulli functions
$$b_q(x) := B_q(\{x\}),$$
where $B_q(x)$ denotes the $q$-th Bernoulli polynomial \cite{apostol50,zagierhallwilson95}.  In this paper, we study the reciprocity phenomenon from a general perspective. Our main objects of study are the \emph{generalized Dedekind sums} defined by
$$S_{\mathbf{f}}(\bm{\nu} | \nu_k) = \sum_{i=1}^{\nu_k-1} \prod_{j\neq k} f_j\left( \frac{i\nu_j }{\nu_k} \right),\quad k=0,\ldots,r,$$
for a vector of distinct nonzero positive integers $\bm{\nu} = (\nu_0,\nu_1,\ldots,\nu_r)$ pairwise coprime, and $\mathbf{f} = (f_0,f_1,\ldots,f_r)\in\cS_1^{r+1}$ where 
$$\cS_{1} := \left\{f:\bR \to \bC : \begin{array}{c}
    \text{$1$-periodic}, \\
    \text{discontinuities at most in } \bZ,\\
    \text{absolutely continuous on }(0,1).
\end{array} \right\}.$$
 
In \Cref{preliminairs}, we review the basic facts on integration, Dedekind sums, and conical zeta values.  In \Cref{recitheoremsec}, we introduce our reciprocity symbol $\cR_{\mathbf{f}}(\bm{\nu})$. Using integration by parts, we obtain an \emph{integral reciprocity formula}. As was observed in \cite{rademacheremil}, this gives \Cref{rademacherreci} in a clean way. The integral reciprocity formula can be reformulated in terms of Fourier analysis. Since every function in $\cS_1$ admits a Fourier series, we derive a \emph{Fourier reciprocity formula}, which connects the theory to polyhedral geometry. We end the section by applying this Fourier reciprocity to the case where the functions $f_j \in \cS_1$ arise from polynomials of degree $q_j\geq 1$ (\Cref{explicitcompsec}). In this setting, reciprocity is closely related to a class of polyhedral zeta values $\zeta_{k,\bm{\nu},\mathbf{q}}$ defined on the lattice $\bm{\nu}_{\bZ}^{\perp}$ of integer vectors orthogonal to $\bm{\nu}$, where $\mathbf{q} = (q_0,q_1, \ldots, q_r)$ is a vector of positive integers. These values provide an explicit reciprocity formula for the Dedekind sums
  $$S_{b_{\mathbf{q}}}(\bm{\nu}|\nu_k) :=   S_{(b_{q_0}\ldots b_{q_r})}(\bm{\nu}|\nu_k)$$
  coming from the periodic Bernoulli functions (\Cref{reciforbq}). In \Cref{zetavaluesandpoly}, we study a two-step decomposition of $\zeta_{k,\bm{\nu},\mathbf{q}}$. First, we establish convergence by comparison with another type of zeta values, denoted $\zeta_{\bm{\nu},\mathbf{q}}$. The key point is that $\zeta_{\bm{\nu},\mathbf{q}}$ converges for all $\mathbf{q}$ and controls the first decomposition. Next, to simplify the situation, we decompose according to the strict integer orthants of $\bR^{r+1}$, i.e., the sets
  $$ \{ (u_0n_0,\ldots,u_rn_{r})  : n_i \in \bN \},\quad (u_0,\ldots,u_r) \in \{-1,1 \}^{r+1},$$
  where $\bN = \{1,2,3,\ldots \}$ is the set of positive integers. This leads to a description in terms of conical zeta values \cite{conicalzeta22,conicalzeta2014, terasoma2004}. This enables us to use desingularization of the cones associated to $\bm{\nu}_{\bZ}^{\perp}$ and express the reciprocity formula in terms of conical zeta values of unimodular cones. In \Cref{dimension2}, we specialize to the case $r=2$, and use the well-known Hirzebruch-Jung algorithm \cite[Ch. 10]{coxlittleschenk2011}, also known as the reversed Euclidean algorithm, to resolve planar cones. This algorithm is used to derive a reciprocity formula for the Dedekind sums $S_{b_{\mathbf{q}}}$ from conical zeta values associated with unimodular cones for each $\mathbf{q}$. 

\section{Preliminaries}\label{preliminairs}

\subsection{Riemann-Stieltjes Integration}
Let $f:[a,b] \to \bC$ be a function, and let $P = \{ t_i \in \bR: a = t_0<t_1<\ldots<t_n =b\}$ be a \emph{partition} of the interval $[a,b]$. The \emph{variation} of $f$ over $P$ is the quantity,
$$V(f,P) := \sum_{i=1}^n |f(t_i)- f(t_{i-1})|.$$
Let $\cP_{a}^b$ be the set of all partitions of $[a,b]$. We say that $f$ is of \emph{bounded variation} on $[a,b]$  if
$$V_{a}^b f = \sup_{P\in \cP_a^b} V(f,P)<\infty.$$
 There is a well-known classification: a function $f$ is of bounded variation if and only if $\re(f)$ and $\im(f)$ can be written as the difference between two increasing functions. As a consequence, functions of bounded variation are differentiable almost everywhere and have at most countably many discontinuities. Thus, we have well-defined two-sided limits,
$$f(x_0^+) := \lim_{x \to x_0^+} f(x),\quad f(x_0^-) := \lim_{x \to x_0^-} f(x).$$ 
The relevance of functions of bounded variation is that they establish a criterion for guaranteeing the existence of the Riemann-Stieltjes integral and the integration by parts formula. Explicitly, set $|P| := \max\{|t_i-t_{i-1}| : i=1,\ldots,n\}$ for any $P\in \cP_a^b$. If $f,g:[a,b] \to \bC$ are functions of bounded variation and have no common discontinuities, then the limits
$$\int_a^b f(x)dg(x) = \lim_{|P|\to 0} \sum_{i=1}^n f(x_i)(g(t_i)-g(t_{i-1})),$$
$$ \int_a^b g(x)df(x) = \lim_{|P|\to 0} \sum_{i=1}^n g(x_i)(f(t_i)-f(t_{i-1})) ,$$
exist for any choice of points $x_i\in [t_{i-i},t_{i}]$ of a partition $P = \{ t_0,\ldots,t_n \} \in \cP_{a}^b$. The following two theorems are well-known (see, for instance, \cite[Ch. 6]{burkill70}).
\begin{theorem}[Integration by Parts]
    \label{integrationbyparts}
    Let $f,g,h:[a,b] \to \bC$ be functions of bounded variation with no common discontinuities. Then,
    $$\int_{a}^b f(x)d(g(x)h(x)) = \int_{a}^b f(x)g(x)dh(x) + \int_{a}^b f(x)h(x)dg(x).$$
\end{theorem}
\qed

As a corollary, for $f(x)=1$ we get the classical integration by parts formula:
$$g(x)h(x)\big{|}_{a}^b = g(b)h(b) - g(a)h(a)= \int_{a}^b g(x)dh(x) + \int_{a}^b h(x)dg(x).$$
We say that $f$ is \emph{absolutely continuous on an interval $I$} if for all $\varepsilon>0$ there exists a $\delta>0$ such that for any collection of disjoint subintervals $(x_1,y_1),\ldots,(x_k,y_k) \subset I$ we have
$$\sum_{j=1}^k |x_j-y_j|<\delta \implies \sum_{j=1}^k |f(x_j)-f(y_j)|<\varepsilon.$$
It is known that an absolutely continuous function is of bounded variation. For a function $f:[a,b]\to \bC$, let $D(f)$ be the set of discontinuities of $f$ in $(a,b)$. Assume $D(f)$ is finite. We say that $f$ is \emph{piecewise absolutely continuous on $I$} if $f$ is absolutely continuous on each connected component of $I\setminus D(f)$.

\begin{theorem}\label{integratorjumps}
    Let $f,g:[a,b] \to \bR$ be functions of bounded variation. If $D(f)\cap D(g) = \emptyset$, $D(f)$ and $D(g)$ are finite, and $g$ is piecewise absolutely continuous on $(a,b)$, then
    $$\int_{a}^b f(x) dg(x) = \int_{a}^b f(x)g'(x) dx + \sum_{x\in D(g)} f(x)(g(x^{+}) - g(x^{-})),$$
    where $g'$ is the derivative of $g$.
\end{theorem}
\qed

\subsection{Dedekind Sums}
Recall that a function $f:\bR \to \bC$ is \emph{periodic} if there exists a real number $\rho > 0$ such that 
$$f(x+\rho) = f(x),\quad \forall x\in \bR.$$
In this case, we say that $f$ is $\rho$-periodic and has \emph{period} $\rho$. A periodic function $f$ is completely determined by its values in $[0,\rho)$. The open interval $(0,\rho)$ is called the \emph{fundamental domain} of $f$. We are interested in the following set of periodic functions defined by
$$\cS_{\rho} := \left\{f:\bR \to \bC : \begin{array}{c}
    \text{$\rho$-periodic}, \\
    \text{discontinuities at most in }\rho\bZ,\\
   \text{absolutely continuous on }(0,\rho).
\end{array} \right\}.$$
  For any function $f\in \cS_{\rho}$ and any positive integer $\nu$ we have the identity
\begin{equation}\label{trivialintegral}
    \int_{0}^{\rho} f(\nu x)dx =  \int_{0}^{\rho} f(x)dx.
\end{equation}
To normalize the situation, we fix all results in the sections that follow for functions $f\in \cS_1$. Indeed, we have a bijection $\cS_1\to \cS_{\rho}$ given by $f(x) \mapsto g(x) := f(\rho^{-1} x)$.

\begin{example}\label{s1examples}The following functions are elements of $\cS_1$.
  \begin{enumerate}
  \item The normalized sine and cosine functions $\mbox{Sin}(x) := \sin(2\pi x)$ and $\mbox{Cos}(x) := \cos(2\pi x) $. In particular, they have no discontinuities. 
      \item The \emph{fractional part} given by $\{ x\} := x-\floor{x}$. Moreover, every function $f\in \cS_1$ relates to the fractional part function by $f(\{ x\}) = f(x).$

\item Let $F:\bR \to \bC$ be any absolutely continuous function. The function given by $f(x) = 
      F(\{x\})$
belongs to $\cS_1$. If $F$ is a polynomial, then we call $f$ a \emph{periodic polynomial function}.
  \end{enumerate}
      \qed
\end{example}
 
\begin{example}\label{examplemodulon}
     For $\nu\in \bN$, each $\cS_{\nu}$ has a \emph{canonical element} given 
    by the \emph{residue mod $\nu$ function} defined as
    $$\{ x\}_{\nu} = x - \nu \floor*{\frac{x}{\nu}},$$
    for any $x\in \bR$. This function has period $\nu$ with discontinuities at $\nu\bZ$, and is determined completely in the interval $[0,\nu)$. An important property is that
\begin{equation}\label{propertymodnu}
        \frac{\{ x\}_{\nu}}{\nu} = \left\{ \frac{x}{\nu}\right\}.
    \end{equation}
    \qed
\end{example}

\begin{definition}
    For $r\geq 1$, let $\mathbf{f} = (f_0,f_1,\ldots,f_r)$ be a vector of functions in $\cS_1$. For any $k \in \{ 0,1,\ldots, r\}$, the \emph{Dedekind sums} associated with this vector are defined as
$$S_{\mathbf{f}}(\bm{\nu} |\nu_k) := \sum_{i=1}^{\nu_k-1} \prod_{j\neq k} f_j\left(\frac{i\nu_j}{\nu_k} \right),$$
for any vector $\bm{\nu} = (\nu_0, \nu_1,\ldots,\nu_r) \in \bN^{r+1}$ of pairwise distinct and pairwise coprime positive integers. Here and throughout,
whenever $k\in\{0,1,\ldots,r\}$, the notation
$
\prod_{j\neq k}
$
means the product over all indices $j\in\{0,1,\ldots,r\}$ with $j\neq k$.
 When $f_j=f$ for all $j$, we simply write $S_{f}$ instead of $ S_{(f,\ldots,f)}$. For $r=0$, we set by convention $S_{f_0}(\bm{\nu} |\nu_0) = \nu_0-1$.
\end{definition}

\begin{remark}
    From \Cref{examplemodulon} we obtain
$$S_{\mathbf{f}}(\bm{\nu} |\nu_k) =S_{\mathbf{f}}(\{ \bm{\nu}\}_{\nu_k}|\nu_k),$$
where, by convention, $\{ \bm{\nu} \}_{\nu_k} = (\{\nu_0\}_{\nu_k},\ldots, 1,\ldots,\{\nu_r\}_{\nu_k})$ with the entry $1$ placed in the $k$-th position.  Thus, for each fixed $\nu_k$, we may restrict ourselves to representatives satisfying $0< \nu_j <\nu_k$ for $j\neq k$.
\end{remark}

\begin{example}
    Important sums of general interest are the Dedekind sums associated with the fractional part
    $$S_{\{x\}}(\bm{\nu} |\nu_k) :=  \sum_{i=1}^{\nu_k-1} \prod_{j\neq k} \left\{\frac{i\nu_j}{\nu_k} \right\}.$$
    Using \Cref{propertymodnu}, these sums satisfy the identity
    $$\nu_k^r S_{\{x\}}(\bm{\nu} |\nu_k) = \sum_{i=1}^{\nu_k-1} \prod_{j\neq k} \left\{i \nu_j\right\}_{\nu_k}.$$
\end{example}

\subsection{Conical Zeta Values}\label{conicalzeta}

   Let $C = \bR_{\geq 0} \mathbf{v}_1 +\ldots+ \bR_{\geq 0} \mathbf{v}_s \subset \bR^{r+1}$ be a polyhedral rational cone,  and let $C^{\circ}:=\mbox{relint}(C)$ denote its relative interior. For any $\mathbf{q} = (q_0,\ldots,q_r) \in \bZ_{\geq 0}^{r+1}$ we define its \emph{conical zeta value} by
    $$\zeta(\mathbf{q},C) := \sum_{\mathbf{n} \in C^{\circ}\cap \bZ^{r+1}} \frac{1}{\mathbf{n}^{\mathbf{q}}},\quad  \frac{1}{\mathbf{n}^{\mathbf{q}}} := \prod_{j=0}^r n_j^{-q_j},$$
     where by convention $0^0 = 1$ and $0^{-q_j}=0$ when $q_{j}\geq 1$. The \emph{reduced conical zeta value} is defined by
$$\zeta^*(\mathbf{q},C) := \sum_{\substack{\mathbf{n} \in C^{\circ}\cap \bZ^{r+1}\\ \gcd(n_0,\ldots,n_r)=1}}  \frac{1}{\mathbf{n}^{\mathbf{q}}}.$$
Denote $|\mathbf{q}|:=q_0+\ldots+q_r$, so both forms are related by
\begin{equation}
    \zeta(\mathbf{q},C)  = \zeta(|\mathbf{q}|)\zeta^*(\mathbf{q},C)
\end{equation}
whenever the series are convergent. 

\begin{remark}\label{closedzeta}
Although our conical zeta values are defined using the relative interior of a cone,
we will also use the corresponding closed contribution of a subdivided cone. More
precisely, for a cone $C$ and subdivision given by a fan $\Sigma$ we set
$$
\zeta_{cl}(\mathbf{q},C,\Sigma)
:=
\sum_{d=1}^{r}\sum_{F \in \Sigma^{(d)}}\zeta(\mathbf{q},F),
$$
where $\Sigma^{(d)}$ denotes the collection of $d$-dimensional cones of $\Sigma$.  Analogously, we can define $\zeta_{cl}^{*}(\mathbf{q},C,\Sigma)$ by running the sum over coprime vectors. If $\Sigma$ is given by the cone itself we just denote $\zeta_{cl}(\mathbf{q},C)$ and $\zeta_{cl}^*(\mathbf{q},C)$ for the previous quantities.
\end{remark}

\begin{example}\label{zetaofvectors}
    For a $1$-dimensional cone $C = \bR_{\geq 0} \mathbf{v}$ generated by a primitive integer vector $\mathbf{v} = (v_0,\ldots ,v_r)$ with $v_j\neq 0$ for all $j$ we get
$$\zeta(\mathbf{q},C) = \sum_{n\in \bN} \prod_{j=0}^r \frac{1}{(v_jn)^{q_j}} = \zeta(|\mathbf{q}|) \prod_{j=0}^r v_j^{-q_j} ,$$
    for $\mathbf{q} \in \bN^{r+1}$. More generally, for $\mathbf{q}\in \bZ_{\geq 0}^{r+1}$, if $v_k=0$ for some  $k=0,\ldots,r$, then by our conventions we have two cases. If $q_k>0$ we have by definition $\zeta(\mathbf{q} , C) = 0$. If $q_k = 0$, then we have the identity
    $$\zeta(\mathbf{q},\bR_{\geq 0} \mathbf{v}) = \zeta(\mathbf{q}^k,\bR_{\geq 0} \mathbf{v}^k),$$
    where $\mathbf{q}^k= (q_0,\ldots,q_{k-1},q_{k+1},\ldots,q_r)$ and $\mathbf{v}^k= (v_0,\ldots,v_{k-1},v_{k+1},\ldots,v_r)$.
\end{example}

\begin{example}
    For the standard positive cone $C = \bR_{\geq 0} \mathbf{e}_0 + \ldots  + \bR_{\geq 0} \mathbf{e}_r$  generated by the canonical basis $\mathbf{e}_0,\ldots ,\mathbf{e}_r$ of $\bR^{r+1}$, we have explicitly
    $$\zeta(\mathbf{q},C) = \sum_{\mathbf{n} \in \bN^{r+1}} \frac{1}{n_0^{q_0}\ldots n_r^{q_r}} = \sum_{n_0\geq 1} \frac{1}{n_0^{q_0}}\ldots \sum_{n_r\geq 1} \frac{1}{n_r^{q_r}} = \prod_{j=0}^r \zeta(q_j),$$
    when $q_j>1$ for all $j$. We can write
    $$\zeta^*(\mathbf{q},C) = \frac{1}{\zeta(|\mathbf{q}|)} \prod_{j=0}^r \zeta(q_j).$$
\end{example}

\section{General Reciprocity Formulas}\label{recitheoremsec}

\subsection{The Integral Reciprocity}

   For a vector of functions $\mathbf{f} = (f_0,f_1,\ldots ,f_r), f_j \in \cS_1$ define, 
$$\Delta f_j:= f_j(1^-) - f_j(0^+),$$
and 
$$\Delta \mathbf{f} := \prod_{j=0}^{r}f_j(1^-) - \prod_{j=0}^{r}f_j(0^+) .$$
The \emph{reciprocity of $\mathbf{f}$} over any $\bm{\nu}  = (\nu_0,\ldots,\nu_r)$ is defined as 
$$\cR_{\mathbf{f}}(\bm{\nu}) := \sum_{k=0}^{r}   S_{\mathbf{f}}(\bm{\nu}|\nu_k)\Delta f_k $$
where we reserve the notation $\cR_{f}(\bm{\nu})$ for the case $\mathbf{f} = (f,\ldots ,f)$. Note that if the functions $f_j$ are continuous for all $j$, then $\cR_{\mathbf{f}}(\bm{\nu})  = 0 $. The following is what we refer to as the \emph{integral reciprocity theorem}.

\begin{theorem}\label{maintheorem}
    Let $\bm{\nu}  = (\nu_0,\ldots,\nu_r)$ be any vector of pairwise distinct and pairwise coprime positive integers. We have,
    $$\cR_{\mathbf{f}}(\bm{\nu}) = -\Delta \mathbf{f}  + \sum_{k=0}^{r}  \nu_k \int_{0}^1 f_k'(\nu_k x)\prod_{j\neq k} f_j(\nu_jx) dx.$$
\end{theorem}
\proof{
Fix a small $\varepsilon>0$, and define
$$I_{\varepsilon} := \int_{\varepsilon}^{1-\varepsilon} d\left( \prod_{j=0}^r f_j(\nu_j x)\right) = \prod_{j=0}^r f_j(\nu_j (1-\varepsilon)) - \prod_{j=0}^r f_j(\nu_j \varepsilon).$$
So, observe that 
$$\lim_{\varepsilon\to 0} I_{\varepsilon} = \prod_{j=0}^{r}f_j(1^-) - \prod_{j=0}^{r}f_j(0^+)  = \Delta \mathbf{f}.$$
On the other hand, since each $f_j$ is of bounded variation in $[\varepsilon,1-\varepsilon]$ and the functions $f_j(\nu_jx)$ have no common discontinuities, we can use integration by parts (\Cref{integrationbyparts}) to get,
$$I_{\varepsilon} = \sum_{k=0}^r \int_{\varepsilon}^{1-\varepsilon} \prod_{j\neq k} f_j(\nu_j x) df_k(\nu_k x).$$
The discontinuities of $f_k(\nu_k x)$ are of the form $i/\nu_k$ for $0\leq i\leq \nu_k$, so $f_k(\nu_k x)$ is piecewise absolutely continuous on $(0,1)$. Using \Cref{integratorjumps}, we obtain,
$$\int_{\varepsilon}^{1-\varepsilon} \prod_{j\neq k} f_j(\nu_j x) df_k(\nu_k x) =  \nu_k\int_{\varepsilon}^{1-\varepsilon}  f_k'(\nu_k x)\prod_{j\neq k} f_j(\nu_j x) dx - \sum_{i=1}^{\nu_k-1} \prod_{j\neq k} f_j\left(\nu_j \frac{i}{\nu_k}\right) (f_k(\nu_k^-) - f_k(\nu_k^+)).$$
Using that $(f_k(\nu_k^-) - f_k(\nu_k^+)) = \Delta f_k$ by periodicity, and taking $\varepsilon \to 0$ the theorem follows.\\
\qed
}

\begin{example}\label{exampleimportant}
    \begin{enumerate} 
        \item The interesting case is when $\Delta f_j = 1$ for all $j$.   For example, the functions $f_j(x) = \{ x\}^{q_j}$ given a collection of integers $q_j\geq 1$. In this case, $f(0^+) = 0$ and $f(1^-) = 1$, so $\Delta \mathbf{f} = 1$ and we have,
        $$\cR_{\mathbf{f}}(\bm{\nu}) = \sum_{k=0}^r S_{\mathbf{f}}(\bm{\nu} | \nu_k) = -1 + \sum_{k=0}^r q_k \nu_k \int_0^1 
 \{\nu_k x\}^{q_k-1} \prod_{j\neq k} \{\nu_j x\}^{q_j}dx.$$

    \item Another important case is when $f_j = b_1$ for all $j=0,\ldots ,r$. We have $f(0^{+}) = -1/2$ and $f(1^{-}) = 1/2$,  $\Delta b_1 = 1$. Also, $b_1'(x) = 1$, so the previous integral reciprocity formula applies as
\begin{equation}
   \cR_{b_1}(\bm{\nu}) =  \sum_{k=0}^{r}   S_{b_1}(\bm{\nu}|\nu_k) = \frac{(-1)^{r+1}-1}{2^{r+1}}  + \sum_{k=0}^{r}  \nu_k \int_{0}^1 \prod_{j\neq k} b_1(\nu_jx) dx.
\end{equation}
Indeed, for $r=2$, we can provide a proof of Dedekind-Rademacher's reciprocity, as observed in \cite[2.D]{rademacheremil}. The key result in the proof is the following identity
$$\int_0^1 b_1(m x) b_1(n x)dx = \frac{\gcd(m,n)}{12mn},\quad \forall m,n \in \bN,$$
due to J. Franel in \cite{franel24}, from which Rademacher's reciprocity follows directly. In his honour, a \emph{Franel integral} is any value of the form
$$\int_{0}^1 b_1(\nu_0x) \ldots  b_1(\nu_rx) dx,$$
for any vector of distinct positive integers $(\nu_0,\ldots ,\nu_r)$. 

\item For a collection $a_0,\ldots ,a_r\in (0,1)$ consider $f_j(x) = \{x \} -a_j$, then we have $\Delta f_j = 1$, and,
$$\Delta \mathbf{f} = \prod_{j=0}^r (1-a_j) - (-1)^{r+1} \prod_{j=0}^r a_j.$$
   For $r=2$, we can obtain a parametric extension of Dedekind-Rademacher's formula. Indeed, using Franel's identity and \Cref{trivialintegral} we get
\begin{align*}
    \int_0^1 (\{\nu_j x\}-a_j)(\{\nu_k x\}-a_k)dx &=\int_0^1 (b_1(\nu_jx)+ 1/2 -a_j)(b_1(\nu_k x) + 1/2-a_k)dx\\
    &= \frac{1}{12\nu_j\nu_k} + (1/2 - a_j)(1/2-a_k). 
\end{align*}
  So, our reciprocity formula is
  \begin{equation}
      \cR_{\mathbf{f}}(\bm{\nu}) =   -1 + \sum_{k=0}^2 a_k -\sum_{j<k} a_ja_k  + \sum_{\substack{j< k\\ l\neq j,k}} \nu_l(1/2 - a_j)(1/2-a_k)+ \frac{\nu_0^2+\nu_1^2+\nu_2^2}{12\nu_0\nu_1\nu_2}.
  \end{equation}
  In particular, when $f_j(x) = f(x) = \{ x\}-a$, for a fixed $0<a<1$ we get
  $$ \cR_{\mathbf{f}}(\bm{\nu})= \frac{(\nu_0^2+\nu_1^2+\nu_2^2)}{12 \nu_0\nu_1\nu_2}  + \left( \frac{1}{2} - a \right)^2(\nu_0+\nu_1+\nu_2) - (1 - 3a + 3a^2).$$
    \end{enumerate}
\end{example}

\subsection{Fourier Approach}\label{fourierrecisec}

 In \cite{berdt2023mcintosh}, the authors use the Fourier expansion of the function $b_1$ to solve a conjecture of McIntosh \cite{96mcintosh} about Franel integrals of order $r=3$. We aim to apply these ideas to study our Dedekind sums. The main reason is that functions of bounded variation always admit a Fourier expansion. Explicitly, for a collection of functions $f_0,\ldots ,f_r\in \cS_1$, we can always assume an expansion
$$\cF{f_j}(x) = \sum_{n\in \bZ} c_{j,n}e(nx),\quad c_{j,n} = c_{n}(f_j) :=  \int_0^1 f_j(x)e(- n x) dx,$$
where $e(x) := \exp(2\pi \iota x)$, and $\iota$ is the imaginary unit with $\iota^2 = -1$. These expansions are characterized by $f_j(x) = \cF{f_j}(x)$ for any $x\in \bR\setminus \bZ$ and are extended to all $\nu\in \bZ$ by
\begin{equation}\label{discmean}
   \cF{f_j}(\nu) = \frac{f_j(\nu^+) + f_j(\nu^-)}{2},
\end{equation}
known as the \emph{principal value} of the series. For any $\mathbf{n}\in \bZ^{r+1}$ denote
$$c_{k,\mathbf{n}} := \prod_{j\neq k } c_{j,n_j}.$$
In what follows, for any pair of vectors $\mathbf{m} = (m_0,\ldots ,m_r),\mathbf{n}= (n_0,\ldots ,n_r)\in \bR^{r+1}$ we have the inner product
$$\bracket{\mathbf{m},\mathbf{n}} = m_0n_0+\ldots +m_rn_r,$$
defining the orthogonal sets
$$\mathbf{m}^{\perp} = \{ \mathbf{n} \in \bR^{r+1} :\bracket{\mathbf{m},\mathbf{n}}= 0 \}.$$
 The following subsets will be relevant,
$$\mathbf{m}^{\perp}_{\Lambda} := \mathbf{m}^{\perp}\cap \Lambda^{r+1},$$
for any discrete subset $\Lambda \subset \bR$. We shall focus on the analysis of the quantities
$$I(\bm{\nu} | \nu_k) := \int_0^1 f_k'(\nu_k x) \prod_{j\neq k} f_j(\nu_j x)dx,$$
 for which the preceding discussion allows us to interchange the Fourier expansions in the principal value sense whenever absolute convergence is not available. 

\begin{corollary}[Fourier Reciprocity]\label{fourierreci} If  $f_k'\in \cS_1$ for all $k$, we have
    $$\cR_{\mathbf{f}}(\bm{\nu}) =   \sum_{k=0}^r \Delta f_k\sum_{\mathbf{n} \in \bm{\nu}^{\perp}_{\bZ}}  \nu_k c_{k,\mathbf{n}} - \Delta \mathbf{f}.$$
\end{corollary}
\proof{
Since $f_k'\in \cS_1$, there exists a Fourier expansion
$$\cF f_k'(x) = \sum_{n\in\bZ} c'_{k,n} e(nx).$$
So, we have explicitly
$$\cF f'_k(\nu_k x)\prod_{j\neq k} \cF f_j(\nu_jx)= \sum_{\mathbf{n} \in \bZ^{r+1}} c'_{k,n_k} \prod_{j\neq k} c_{j,n_j} e\left( \bracket{\bm{\nu},\mathbf{n}}x \right) ,$$
converging point-wise for any $x$ and uniformly on any compact subset of  $\bR\setminus \bZ$, so integrals can be interchanged with all sums understood in the above principal value sense. Recall that 
    $$ \int_{0}^1 e\left( \bracket{\bm{\nu},\mathbf{n}}x \right) dx=  \left\{ \begin{matrix}
        1, & \bracket{\bm{\nu},\mathbf{n}}=0,\\
        0, &  \bracket{\bm{\nu},\mathbf{n}} \neq 0,
    \end{matrix}\right.$$
so we can compute directly,
$$I(\bm{\nu} | \nu_k) = \sum_{\mathbf{n} \in \bm{\nu}^{\perp}_{\bZ}} c'_{k,n_k} \prod_{j\neq k} c_{j,n_j}. $$
Therefore, integral reciprocity can be rewritten as follows,
\begin{equation}\label{recioff}
 \cR_{\mathbf{f}}(\bm{\nu}) = \sum_{k=0}^{r} \Delta f_k S_{\mathbf{f}}(\bm{\nu} |\nu_k)  = -\Delta \mathbf{f} + I(\nu),
\end{equation}
where
$$I(\bm{\nu} ):=  \sum_{k=0}^r \nu_k I(\bm{\nu} |\nu_k) = \sum_{k=0}^r \sum_{\mathbf{n} \in \bm{\nu}_{\bZ}^{\perp}} \nu_k c'_{k,n_k} \prod_{j\neq k} c_{j,n_j}.$$
Observe that $c_{k,0}' = \Delta f_k$, and for any $n\neq 0$
\begin{align*}
    c_{k,n} = \int_{0}^1 f_k(x)e(-n x)dx  &= \left.-\frac{f_k(x)e(-n x)}{2\pi \iota n}\right|_{x=0}^{x=1} + \frac{1}{2\pi \iota n} \int_0^1 f_k'(x) e(- n x)dx \\
    &= \frac{1}{2\pi \iota n} \left(c'_{k,n} - \Delta f_k \right).
\end{align*}
So $c_{k,n}' = (2\pi \iota n) c_{k,n} + \Delta f_k$ for all $n$, and now we can rewrite,
$$I(\bm{\nu} |\nu_k)  = \Delta f_k \sum_{\mathbf{n} \in \bm{\nu}^{\perp}_{\bZ}}  \prod_{j\neq k} c_{j,n_j} +  2\pi\iota  \sum_{\substack{\mathbf{n} \in \bm{\nu}^{\perp}_{\bZ}}} n_k  \prod_{j = 0}^r c_{j,n_j}.$$
Summing over $k=0,\ldots ,r$, observe that in terms of principal values we have,
$$ 2\pi\iota \sum_{k=0}^r \sum_{\substack{\mathbf{n} \in \bm{\nu}^{\perp}_{\bZ} }} \nu_k n_k  \prod_{j = 0}^r c_{n_j} = 2\pi\iota  \sum_{\substack{\mathbf{n} \in \bm{\nu}^{\perp}_{\bZ}}} \prod_{j = 0}^r c_{n_j} \sum_{k=0}^r \nu_k n_k   = 0.$$ 
So, we get the identity
\begin{equation}\label{Ivtupendo}
    I(\bm{\nu} ) = \sum_{k=0}^r \Delta f_k\sum_{\mathbf{n} \in \bm{\nu}^{\perp}_{\bZ}}   \nu_k c_{k,\mathbf{n}},
\end{equation}
from which the result follows.\\
\qed}

\begin{remark}\label{remarksimportant}
    Important particular cases of the previous result are the following.
    \begin{enumerate}
        \item Assume that some of the functions $f_j$ in $\mathbf{f}$ are continuous in $\bR$, i.e., $\Delta f_j = 0$. Let us denote $\cC = \{j : f_j \text{ is continuous} \}$, so reciprocity looks like
        $$\cR_{\mathbf{f}}(\bm{\nu}) = \sum_{k \not\in \cC} S_{\mathbf{f}}(\bm{\nu}|\nu_k)\Delta f_k   = \sum_{k \not\in \cC} \sum_{\mathbf{n} \in \bm{\nu}^{\perp}_{\bZ}}  \nu_k \Delta f_k c_{k,\mathbf{n}}  -  \prod_{k\in \cC} f_k(0) \Delta \tilde{\mathbf{f}},$$
        where $\tilde{\mathbf{f}}$ is the vector of discontinuous functions of $\mathbf{f}$. In particular, if $f_k(x) = \{x\}-a$ for any $a\in \bR$ and $f_j$ continuous for all $j \neq k$, then
    $$\cR_{\mathbf{f}}(\bm{\nu}) = S_{\mathbf{f}}(\bm{\nu}|\nu_k)  =  \nu_k \sum_{\mathbf{n} \in \bm{\nu}^{\perp}_{\bZ}}   c_{k,\mathbf{n}}  -  \prod_{j\neq k } f_j(0).$$
    This formula is especially simple for standard Fourier test functions such as
$e(x) = \exp(2\pi \iota x)$, $\mbox{Cos}(x) = \cos(2\pi x)$, and $\mbox{Sin}(x) = \sin(2\pi x)$, since their Fourier
series have finite support. For instance, taking \(f_j(x)=e(x)\) for all \(j\neq k\),
one obtains
    $$
S_{\mathbf f}(\nu\mid \nu_k)
=
\begin{cases}
\nu_k-1, & \text{if } \nu_k\mid \sum_{j\neq k}\nu_j,\\
-1, & \text{otherwise}.
\end{cases}
$$
The cosine and sine cases are obtained similarly by decomposing them into exponentials. Indeed,
$$ S_{\mathbf{f}}(\bm{\nu}|\nu_k) = \left\{ \begin{matrix}
    \dfrac{\nu_k}{2^r}|U_k|  -1 & \text{ if $f_j(x) = \mbox{Cos}(x)$ for all \(j\neq k\)},\\[1em]
    \dfrac{\nu_k}{(2\iota)^{r}}    \sum_{\mathbf{u}\in U_k} \mbox{sgn}(\mathbf{u}) & \text{ if $f_j(x) = \mbox{Sin}(x)$ for all \(j\neq k\),}
\end{matrix} \right. $$
where $$U_k = \left\{\mathbf{u} \in \{-1,1\}^r :  \nu_k \big| \bracket{\mathbf{u}, \bm{\nu}^{k}} \right\},$$
with $\bm{\nu}^{k} = (\nu_0,\ldots ,\nu_{k-1},\nu_{k+1},\ldots ,\nu_r)$, and $\mbox{sgn}(\mathbf{u})$ is the product of its components.

    \item For any $f\in \cS_1$ not constant with $f'\in \cS_1$ and $f(0^+) \neq f(1^-)$, we have
    $$\sum_{k=0}^{r} S_f(\bm{\nu}|\nu_k) =   \sum_{k=0}^r \sum_{\mathbf{n} \in \bm{\nu}^{\perp}_{\bZ}}  \nu_k c_{k, \mathbf{n}} - \frac{\Delta \mathbf{f}}{\Delta f} ,$$
    where explicitly
    $$\frac{\Delta \mathbf{f}}{\Delta f} =\sum_{k=0}^rf(1^-)^{k}f(0^+)^{r-k}.$$
    \end{enumerate}
\end{remark}

  \subsection{Reciprocity for Polynomial Dedekind Sums}\label{explicitcompsec}

For polynomials $F_0,\ldots ,F_r \in \bC[x]$ of degrees $q_j\geq 1$, consider their associated periodic functions  $f_j\in \cS_1$ defined by $f_j(x) : = F_j(\{x\})$. Since the study of the Dedekind sums $S_{\mathbf{f}}(\bm{\nu}|\nu_k)$ reduces to the study of such sums for a basis of $\bC[x]$, we use the Bernoulli polynomials $B_q(x)$ defined by the generating function
$$\frac{te^{xt}}{e^t - 1} = \sum_{q \geq 0} B_q(x)\frac{t^q}{q!}.$$
In $\cS_1$, we have the \emph{periodic $q$-th Bernoulli polynomials} $b_q(x) = B_q(\{x \})$, and their Fourier coefficients are computed as
\begin{equation}\label{FCqbern}
    c_n(b_q) = \int_0^1 B_q(x) e(-nx) dx =  \left\{\begin{array}{cc}
      c_0(b_0) = 1,\quad  c_n(b_0) = 0,\, n\neq 0,& \\[3mm]
      c_0(b_q) = 0, \quad c_n(b_q)=-\dfrac{q!}{(2\pi \iota n)^q}, &  n\neq 0,\, q>0.
 \end{array} \right.
\end{equation}
The relation with the power functions is
$$x^q = \frac{1}{q+1}\sum_{s=0}^q \binom{q+1}{s} B_{s}(x).$$
For any $\mathbf{q} = (q_0, \ldots, q_r) \in (\bZ_{\geq 0})^{r+1}$ we will denote $S_{\mathbf{q}}$ and $S_{b_{\mathbf{q}}}$ for the Dedekind sums associated to $S_{(\{x\}^{q_0},\ldots, \{x\}^{q_r})}$ and $S_{(b_{q_0},\ldots, b_{q_r})}$ respectively. 
\begin{remark}
If one or more components of \(\mathbf q\) vanish, the associated Dedekind
sums reduce to the corresponding sums attached to the subvector obtained
by deleting the zero components. These reduced cases should be kept in mind
in what follows.
\end{remark}
From the previous relation, it is not too difficult to check that
$$\Delta \{x\}^{q_k}S_{\mathbf{q}}(\bm{\nu}|\nu_k) = \prod_{j=0}^{r} \frac{1}{q_j+1} \sum_{\substack{\mathbf{s}\in \bZ^{r+1}\\ 0\leq s_j\leq q_j}} A_{\mathbf{s}} \Delta b_{s_k} S_{b_{\mathbf{s}}}(\bm{\nu}|\nu_k),\quad A_{\mathbf{s}} := \prod_{j=0}^r \binom{q_j+1}{s_j}.$$
In this way, if $\cR_{\mathbf{q}}(\bm{\nu})$ denotes the reciprocity of the sums $S_{\mathbf{q}}$, then in terms of reciprocity we have,
\begin{equation}\label{xqrecibq}
    \cR_{\mathbf{q}}(\bm{\nu}) = \prod_{j=0}^{r} \frac{1}{q_j+1} \sum_{\substack{\mathbf{s}\in \bZ^{r+1}\\ 0\leq s_j\leq q_j}} A_{\mathbf{s}} \cR_{b_{\mathbf{s}}}(\bm{\nu}).
\end{equation}
 Since the periodic Bernoulli functions are continuous for $q>1$, we have
$$\Delta b_{q} = \left\{\begin{matrix}
    1 & q=1,\\
    0 & q\neq 1,
\end{matrix}\right. $$
and thus,
$$  \cR_{b_{\mathbf{q}}}(\bm{\nu}) =  \sum_{ q_k = 1} S_{b_{\mathbf{q}}}(\bm{\nu}|\nu_k).$$
In particular, this reduces the index set of the right-hand side sum in \Cref{xqrecibq} into the set $\{\mathbf{s} \in \bZ^{r+1} : 0\leq s_j \leq q_j, \exists s_j = 1 \}$. We adopt the notation
$$\mathbf{q}_k = (q_0,\ldots ,q_{k-1},0,q_{k+1},\ldots ,q_r),\quad  I_{\mathbf{q},k} :=  \frac{(-1)^r}{(2\pi \iota)^{|\mathbf{q}_k|}}\prod_{j\neq k} q_j!,$$
$$\mathbf{1}_{\mathbf{q}} := \#\{q_j \in \mathbf{q} : q_j = 1\},\quad |\overline{\mathbf{q}}| := |\mathbf{q}|-1,\quad I_{\overline{\mathbf{q}}} := \frac{(-1)^r}{(2\pi \iota)^{|\overline{\mathbf{q}}|}}\prod_{q_j>1} q_j!.$$
Moreover, observe that $I_{\mathbf{q},k} = I_{\overline{\mathbf{q}}}$ for all $k$ such that $q_k = 1$. Now, assume that all components of \(\mathbf q\) are positive. From \Cref{FCqbern}, for any $k=0,\ldots ,r$ and any $\mathbf{n}\in\bZ^{r+1}$ we can
explicitly compute the coefficients $c_{k,\mathbf{n}}$ for the vector of functions $\mathbf{f} = (b_{q_0},\ldots,b_{q_{r}})$, namely
 $$c_{k,\mathbf{n}} = \prod_{j\neq k} c_{n_j}(b_{q_j}) =I_{\mathbf{q},k}\frac{1}{\mathbf{n}^{\mathbf{q}_k}}.$$
This motivates the definition
$$\zeta_{k,\bm{\nu},\mathbf{q}} := \sum_{\mathbf{n} \in \bm{\nu}^{\perp}_{\bZ} } \frac{1}{\mathbf{n}^{\mathbf{q}_k}}.$$
In the present setting, this quantity will be referred to as a  \emph{multiple zeta value}.
\begin{remark}
The term \emph{multiple zeta value} is used here in the specific sense
for the sums defined above. This should be
distinguished from the classical multiple zeta values studied by Hoffman \cite{hoffmanmzv92},
Zagier \cite{multzagier12}, and others, although the terminology is meant to emphasize the
analogous zeta-type nature of these sums.
\end{remark}
From (1) in \Cref{remarksimportant}, we summarize the previous discussion into the following.

\begin{corollary}[Reciprocity for Bernoulli Polynomials]\label{reciforbern}
         For any positive integer vector $\mathbf{q}$, we have\begin{equation}\label{reciforbq}
    \cR_{b_{\mathbf{q}}}(\bm{\nu}) =  I_{\overline{\mathbf{q}}} \sum_{q_k = 1}\nu_k \zeta_{k,\bm{\nu},\mathbf{q}} -  \frac{1 -  (-1)^{\mathbf{1}_{\mathbf{q}}}}{2^{\mathbf{1}_{\mathbf{q}}}}  \prod_{q_j\neq 1} B_{q_j}.
\end{equation}
\end{corollary}
\qed

Since the sums $\zeta_{k,\bm{\nu},\mathbf{q}}$ come from a Fourier analysis approach, if they are not absolutely convergent, then we consider their principal values. So, if $|\mathbf{q}_k|$ is odd, then the antipodal map $\mathbf{n}\mapsto -\mathbf{n}$ gives,
$$\zeta_{k,\bm{\nu},\mathbf{q}} = \sum_{\mathbf{n}\in \bm{\nu}^{\perp}_{\bZ}\cap \{n_k\geq0\} } \frac{1}{\mathbf{n}^{\mathbf{q}_k}} + (-1)^{|\mathbf{q}_k|} \sum_{\mathbf{n} \in  \bm{\nu}^{\perp}_{\bZ}\cap \{n_k\geq0\}} \frac{1}{\mathbf{n}^{\mathbf{q}_k}} = 0.$$
The condition $|\mathbf{q}_k|$ being odd for all $k=0,\ldots ,r$ implies that $r$ and all the $q_j$ are odd. And since $B_q = 0$ for all odd integers $q>2$, as a direct consequence, the following holds.
    \begin{corollary} 
   If $\mathbf{q} = (q_0,\ldots ,q_r)$ satisfies $q_k>1$ for all $k$ or  $|\mathbf{q}_k|$ is odd for all $k=0,\ldots ,r$, then $\cR_{b_{\mathbf{q}}}(\bm{\nu}) = 0.$
\end{corollary}
\qed

\section{Zeta Values and Polyhedra}\label{zetavaluesandpoly}

From the previous section, the reciprocity $\cR_{b_{\mathbf{q}}}(\bm{\nu})$ for $r=1$ is nontrivial when we assume $\mathbf{q} = (1,q)$ where $q$ is a positive integer. Then, we have explicitly
$$\cR_{b_{\mathbf{q}}}(\nu_0,\nu_1) = -\frac{\nu_0 q!}{(2\pi \iota)^q} \sum_{\mathbf{n} \in \nu^{\perp}_{\bZ}} \frac{1}{n_1^q} - \delta(q)\frac{\nu_1}{(2\pi \iota)} \sum_{\mathbf{n} \in \nu^{\perp}_{\bZ}} \frac{1}{n_0} - B(q) ,  $$
where 
$$\delta(q) = \left\{\begin{matrix}
    1 & q=1,\\
    0 & q>1,
\end{matrix}\right. \quad B(q) = \left\{\begin{matrix}
    0 & q=1,\\
    B_q & q>1.
\end{matrix}\right.$$
The set
$$\nu^{\perp}_{\bZ} = \{(n_0,n_1) : \nu_0n_0 + \nu_1n_1 = 0 \},$$
has well-known solutions given by $(n_0,n_1) = n(\nu_1,-\nu_0)$ for all $n\in \bZ$. So,  reciprocity becomes
$$\cR_{b_{\mathbf{q}}}(\nu_0,\nu_1) = -\frac{(-1)^q q!}{\nu_0^{q-1}(2\pi \iota)^q}  \sum_{\substack{n\in \bZ \\ n\neq 0}} \frac{1}{n^q} - \delta(q)\frac{1}{(2\pi \iota)} \sum_{\substack{n\in \bZ \\ n\neq 0}}  \frac{1}{n} - B(q).$$
Recall that
$$\sum_{\substack{n\in \bZ \\ n\neq 0}} \frac{1}{n^q} = \left\{\begin{matrix}
    0 & q \text{ odd},\\
    2\zeta(q) & q \text{ even},
\end{matrix}\right. $$
where
$$\zeta(q) := \sum_{n\in\bN} \frac{1}{n^q}$$
is the Riemann zeta function. The remarkable identity
$$\zeta(q) = (-1)^{\frac{q}{2}+1}\frac{B_{q}(2\pi)^{q}}{2q!},$$
holds for all $q$ even. Therefore, reciprocity for $r=1$ is
\begin{equation}
    \cR_{b_{\mathbf{q}}}(\nu_0,\nu_1) = \left\{\begin{matrix}
    - B(q) & q \text{ odd},\\[5pt]
    \dfrac{(-1)^{1 + \frac{q}{2}}q!}{\nu_0^{q-1}2^{q-1}\pi^q}  \zeta(q)- B(q) & q \text{ even}
\end{matrix}  \right.=\left\{\begin{matrix}
    - B(q) & q \text{ odd},\\[5pt]
    \left(\nu_0^{1-q}- 1\right) B_q & q \text{ even}.
\end{matrix}  \right.
\end{equation}
In what follows, we treat the case $r\geq 2$ for the multiple zeta values $\zeta_{k,\bm{\nu},\mathbf{q}}$. This will lead us into a reciprocity formula for $r=2$ in the next section.

\subsection{The Multiple Zeta Values.}

 In what follows,  $\bZ^{\times}$ will denote the set of nonzero integers. Define
$$\bm{\nu}^{k,\perp}_{\bZ} := \bm{\nu}^{\perp}_{\bZ}\cap \{ n_k=0\} \cong (\bm{\nu}^{k})^{\perp}_{\bZ},$$
where $\bm{\nu}^{k} = (\nu_0,\ldots ,\nu_{k-1},\nu_{k+1},\ldots ,\nu_r)$. Thus, we can decompose
    $ \zeta_{k,\bm{\nu},\mathbf{q}} = Y_{k,\bm{\nu},\mathbf{q}} +  Z_{k,\bm{\nu},\mathbf{q}},$ where
    $$Y_{k,\bm{\nu},\mathbf{q}} := \sum_{\mathbf{n}\in\bm{\nu}^{\perp}_{\bZ^{\times}}} \frac{1}{\mathbf{n}^{\mathbf{q}_k}},\quad Z_{k,\bm{\nu},\mathbf{q}} := \sum_{\mathbf{n}\in\bm{\nu}^{k,\perp}_{\bZ}} \frac{1}{\mathbf{n}^{\mathbf{q}_k}}.$$
Observe that, in both cases, the corresponding principal values vanish when $|\mathbf{q}_k|$ is odd. Hence, these terms do not contribute to the
evaluations. However, the sums $Y_{k,\bm{\nu},\mathbf{q}}$ do not converge on their own if $q_j = 1$ for some $j\neq k$. The right side of \Cref{reciforbern} tells us how these sums should be combined to obtain
convergent expressions. For this purpose, define the following multiple zeta value
$$\zeta_{\bm{\nu}, \mathbf{q}} := \sum_{\mathbf{n}\in\bm{\nu}^{\perp}_{\bZ^{\times}}} \frac{1}{\mathbf{n}^{\mathbf{q}}}.$$

\begin{theorem}\label{convergence}
    The sums $\zeta_{\bm{\nu}, \mathbf{q}}$ converge absolutely for all $\mathbf{q}\in \bN^{r+1}$. Furthermore,
    $$Z_{k,\bm{\nu}, \mathbf{q}} = \zeta_{\bm{\nu}^k,\mathbf{q}^k},$$
    $$\sum_{q_k = 1} \nu_k Y_{k,\bm{\nu}, \mathbf{q}} = -\sum_{q_k>1} \nu_k \zeta_{\bm{\nu}, \mathbf{q} - \mathbf{e}_k}.  $$
\end{theorem}
\proof{First observe that for any $\mathbf{q}\in \bN^{r+1}$, we have
$$|\zeta_{\bm{\nu}, \mathbf{q}}| \leq  \sum_{\mathbf{n}\in\bm{\nu}^{\perp}_{\bZ^{\times}}}  \left| \frac{1}{\mathbf{n}^{\mathbf{q}}}\right| \leq\sum_{\mathbf{n}\in\bm{\nu}^{\perp}_{\bZ^{\times}}}  \left| \frac{1}{\mathbf{n}^{(1,\ldots,1)}}\right|.$$
For any $\mathbf{n} = (n_0,n_2,\ldots,n_r)$, denote $M_1 = \sum_{j\neq 1} \nu_j n_j$, so explicitly, we have
\begin{align*}
\sum_{\mathbf{n}\in\bm{\nu}^{\perp}_{\bZ^{\times}}}  \left| \frac{1}{\mathbf{n}^{(1,\ldots,1)}}\right| &= \sum_{n_0 \in \bZ^{\times}} \frac{1}{|n_0|} \sum_{\substack{(n_j)_{j\neq 0} \in (\bZ^{\times})^{r} \\ \sum_{j\neq 0} \nu_j n_j = -\nu_0n_0}} \prod_{j\neq 0} \frac{1}{|n_j|}\\
    &= \nu_1\sum_{n_0 \in \bZ^{\times}} \frac{1}{|n_0|} \sum_{\substack{(n_j)_{j\neq 0,1} \in (\bZ^{\times})^{r-1} \\ M_1 \neq 0,\, \nu_1|M_1}} \frac{1}{|\sum_{j\neq 1} \nu_j n_j |} \prod_{j\neq 0,1} \frac{1}{|n_j|},
\end{align*}
Now, for any $\mathbf{n} = (n_0,n_3,\ldots,n_r)$ denote $M_{12} =  \sum_{j\neq 1,2} \nu_j n_j$ and distinguish between three cases:
\begin{enumerate}
    \item If $|\nu_{2}n_{2}| \leq |M_{12}|/2$, then 
    $$\left|\sum_{j\neq 1} \nu_j n_j\right| = |\nu_{2}n_{2} + M_{12}| \geq \left|M_{12}\right| - |\nu_{2}n_{2}| \geq \frac{|M_{12}|}{2}.$$
    So, we have
    $$\sum_{\substack{n_{2} \neq 0\\ |\nu_{2}n_{2}| \leq |M_{12}|/2}} \frac{1}{|\sum_{j\neq 1} \nu_j n_j | |n_{2}|} \leq  \frac{2}{|M_{12}|} \sum_{\substack{n_{2} \neq 0\\ |\nu_{2}n_{2}| \leq |M_{12}|/2}} \frac{1}{|n_{2}|} =  
       O(\log(|M_{12}|)|M_{12}|^{-1}).$$
   
    \item If $|M_{12}|/2<|\nu_{2}n_{2}|<2|M_{12}|$, then 
    $$||\nu_{2}n_{2}| - |M_{12}|| \leq |M_{12}|.$$
    So, we have
    \begin{align*}
        \sum_{\substack{n_{2} \neq 0\\ |M_{12}|/2 < |\nu_{2}n_{2}| < 2|M_{12}|}} \frac{1}{|\sum_{j\neq 1} \nu_j n_j | |n_{2}|}  &\leq  \frac{2\nu_2}{|M_{12}|} \sum_{\substack{n_{2} \neq 0\\ |M_{12}|/2 < |\nu_{2}n_{2}| < 2|M_{12}|}} \frac{1}{|\nu_{2}n_2 + M_{12}|}\\
        &\leq  \frac{2\nu_2}{|M_{12}|} \sum_{\substack{n_{2} \neq 0\\ 0 < |k| < |M_{12}|}} \frac{1}{|k|}\\
        &= 
       O(\log(|M_{12}|)|M_{12}|^{-1}).
    \end{align*}
    
    \item If $|\nu_{2}n_{2}| \geq 2|M_{12}|$, then
    $$\left|\sum_{j\neq 1} \nu_j n_j\right| = |\nu_{2}n_{2} + M_{12}| \geq |\nu_{2}n_{2}|  - \left|M_{12}\right|  \geq \frac{|\nu_{2}n_2|}{2}. $$
    So, we have
    $$\sum_{\substack{n_{2} \neq 0\\ |\nu_{2}n_{2}| \geq 2|M_{12}|}} \frac{1}{|\sum_{j\neq 1} \nu_j n_j | |n_{2}|} \leq \frac{2}{\nu_2} \sum_{\substack{n_{2} \neq 0\\ |\nu_{2}n_{2}| \geq 2|M_{12}|}} \frac{1}{|n_{2}|^2} = 
       O(|M_{12}|^{-1}) .$$
\end{enumerate}
Observe that the contribution of the terms with $M_{12}=0$ is lower-dimensional and converges by induction. Therefore, the proof of the convergence of $|\zeta_{\bm{\nu}, \mathbf{q}}|$ is reduced to checking the convergence of
$$\sum_{n_0 \in \bZ^{\times}} \frac{1}{|n_0|} \sum_{\substack{(n_j)_{j\neq 0,1,2} \in (\bZ^{\times})^{r-2} \\ M_{12} \neq 0}} \frac{1}{|M_{12} |^{1-\varepsilon}} \prod_{j\neq 0,1,2} \frac{1}{|n_j|},$$
for every sufficiently small $\varepsilon> 0$ by the fact that we can always compare $\log(x) = O(x^{\varepsilon})$. This process can be done analogously for $M_{123} = \sum_{j\neq 1,2,3} \nu_jn_j$ and so on until $M_{1\ldots r} = \nu_0n_0$. Therefore, we reduce the convergence of $|\zeta_{\bm{\nu}, \mathbf{q}}|$ to the convergence of
$$\sum_{n_0 \in \bZ^{\times}} \frac{1}{|n_0||\nu_0 n_0|^{1-\varepsilon}} = \frac{2}{\nu_0^{1-\varepsilon}} \zeta(2 -\varepsilon),$$
assuring the convergence for small $\varepsilon>0$. The identity for $Z_{k,\bm{\nu},\mathbf{q}}$ is direct by definition. For the last identity, without loss of generality, let us assume that $q_0,q_1,\ldots,q_s>1$ and $q_k=1$ for $k>s$. Now we can write,
$$\sum_{q_k = 1} \nu_k   Y_{k,\bm{\nu},\mathbf{q}} = \sum_{\mathbf{n} \in \bm{\nu}^{\perp}_{\bZ^{\times}} } \sum_{k = s+1 }^r \frac{\nu_k}{\mathbf{n}^{\mathbf{q}_k}} =  \sum_{\mathbf{n} \in \bm{\nu}^{\perp}_{\bZ^{\times}} } \frac{1}{\mathbf{n}^{\mathbf{q}}}  \sum_{k = s+1 }^r \nu_kn_k.$$
Since $\bracket{\bm{\nu},\mathbf{n}} = 0$, the sum is modified into
$$\sum_{q_k = 1} \nu_k   Y_{k,\bm{\nu},\mathbf{q}} = -\sum_{\mathbf{n} \in \bm{\nu}^{\perp}_{\bZ^{\times}} } \frac{1}{\mathbf{n}^{\mathbf{q}}}  \sum_{k = 0}^s \nu_kn_k = -\sum_{q_k>1} \nu_k \sum_{\mathbf{n} \in \bm{\nu}^{\perp}_{\bZ^{\times}} } \frac{n_k}{\mathbf{n}^{\mathbf{q}}}  = - \sum_{q_k>1}\nu_k \zeta_{\bm{\nu}, \mathbf{q} - \mathbf{e}_k} .$$
\qed
}

\begin{remark}\label{vanishremark}
    \begin{enumerate}

     \item The previous discussion says that  in the case $\mathbf{q} = (1,\ldots ,1)$ we must have
    $$\sum_{k=0}^r\nu_k  Y_{k,\bm{\nu},\mathbf{q}} = 0.$$
    Thus, we get explicitly
\begin{equation}
        \sum_{k=0}^r \nu_k \zeta_{k,\bm{\nu},\mathbf{q}}= \sum_{k=0}^r\nu_k Z_{k,\bm{\nu}, \mathbf{q}}.
    \end{equation}

        \item For all $\mathbf{q}\in \bN^{r+1}$ we have
    $$Z_{k,\bm{\nu},\mathbf{q}} =  (-1)^{r}\prod_{j\neq k}\frac{(2\pi \iota)^{q_j}}{q_j!} I(k,\bm{\nu} ,\mathbf{q}),$$
    where
    $$I(k, \bm{\nu}, \mathbf{q}) := \int_0^1 \prod_{j\neq k}b_{q_j}\left(\nu_jx\right)dx.$$
    Moreover, $Z_{k,\bm{\nu},\mathbf{q}}$ is a rational multiple of $(2\pi)^{|\mathbf{q}_k|}$. This follows directly from the proof of \Cref{fourierreci} when we assume $q_k = 1$. Since the $b_{q_j}$ are defined from rational polynomials, the integrals $I(k, \bm{\nu}, \mathbf{q})$ must be rational numbers.

    \end{enumerate}
\end{remark}

The following reduction will be useful for explicit computations in the next section. Let us translate the situation into the first orthant $\bN^{r+1}$. For each $\mathbf{u} = (u_0,u_1,\ldots ,u_r) \in \{ -1,1\}^{r+1}$ the Hadamard product with $\bm{\nu}$ is just $\mathbf{u}\bm{\nu} := (u_0\nu_0,\ldots ,u_r\nu_r).$ We have a bijection,
$$\bigsqcup_{\mathbf{u}\in \{-1,1 \}^{r+1}}(\mathbf{u}\bm{\nu})^{\perp}_{\bN} \to \bm{\nu}_{\bZ^{\times}}^{\perp},\quad (n_0,\ldots ,n_r) \in  (\mathbf{u}\bm{\nu})^{\perp}_{\bN} \mapsto \mathbf{u}\mathbf{n} = (u_0n_0,\ldots ,u_rn_r) .$$
Thus, we get
$$ Y_{k,\bm{\nu},\mathbf{q}}= \sum_{\mathbf{u} \in \{ -1,1\}^{r+1}} \sigma_k(\mathbf{u}) Y_{k,\bm{\nu},\mathbf{q},\mathbf{u}},\quad Y_{k,\bm{\nu},\mathbf{q},\mathbf{u}} = \sum_{\mathbf{n} \in (\mathbf{u}\bm{\nu})^{\perp}_{\bN} } \frac{1}{\mathbf{n}^{\mathbf{q}_k}},$$
where $\sigma_k(\mathbf{u}) = \prod_{j\neq k} u_j^{q_j}$. The same can be done for $Z_{k,\bm{\nu},\mathbf{q}}$ using the identification $\bm{\nu}_{\bZ}^{k,\perp} \cong (\bm{\nu}^{k})^{\perp}_{\bZ^{\times}}$, so we can rewrite,
$$ Z_{k,\bm{\nu},\mathbf{q}}= \frac{1}{2}\sum_{\mathbf{u} \in \{ -1,1\}^{r+1}} \sigma_k(\mathbf{u})Z_{k,\bm{\nu},\mathbf{q},\mathbf{u}},\quad Z_{k,\bm{\nu},\mathbf{q},\mathbf{u}} = \sum_{\mathbf{n} \in (\mathbf{u}\bm{\nu}^{k})^{\perp}_{\bN}} \frac{1}{\mathbf{n}^{\mathbf{q}_k}},$$
where the $2^{-1}$ comes from avoiding repeated sums since the variable $n_k$ is avoided. Indeed, it is possible to remove the $2^{-1}$ factor by redefining the sum over the $\mathbf{u} \in\{-1, 1\}^r$. As above, the cases $q_j = 1$ make the sums $Y_{k,\bm{\nu},\mathbf{q},\mathbf{u}}$ diverge. However, we can combine them as in \Cref{convergence} to establish convergence. Therefore, in what follows, we treat them as a single value without regard to convergence. When $q_{j}>1$ for all $j$, we get the following estimates.

\begin{prop}
If $\nu_l = \min\{\nu_j : j\neq k \}$, then the following holds,
$$Z_{k,\bm{\nu},\mathbf{q},\mathbf{u}} \leq \prod_{j\neq k} \zeta(q_j),\quad  Y_{k,\bm{\nu},\mathbf{q},\mathbf{u}} \leq 2 \nu_l^{q_l} \prod_{j\neq k} \zeta(q_j),$$
for all $\mathbf{q} \in \bN_{>1}^{r+1}$. 
\end{prop}
\proof{ We have
$$Z_{k,\bm{\nu},\mathbf{q},\mathbf{u}}  = \sum_{\mathbf{n} \in (\mathbf{u}\bm{\nu}^{k})^{\perp}_{\bN}} \prod_{j\neq k} \frac{1}{n_j^{q_j}} \leq \sum_{(n_j)_{j\neq k} \in \bN^{r}} \prod_{j\neq k} \frac{1}{n_j^{q_j}} = \prod_{j\neq k}\zeta(q_j).$$
For $\mathbf{n} = (n_0,\ldots ,n_r)\in  (\mathbf{u}\bm{\nu})^{\perp}_{\bN}$ choose an $l\neq k$ and write $u_l\nu_ln_l = - \sum_{j\neq l} u_j\nu_j n_j$. Then, we get
$$Y_{k,\bm{\nu},\mathbf{q},\mathbf{u}} =\sum_{\mathbf{n}\in  (\mathbf{u}\bm{\nu})^{\perp}_{\bN} }  \frac{1}{\mathbf{n}^{\mathbf{q}_k}} =  (-u_l\nu_l)^{q_l}\sum_{\mathbf{n}\in  (\mathbf{u}\bm{\nu})^{\perp}_{\bN} } \left( \frac{1}{  \sum_{j\neq l} u_j\nu_j n_j} \right)^{q_l}\prod_{j\neq k,l} \frac{1}{n_j^{q_j}}.$$
Then, we have
$$Y_{k,\bm{\nu},\mathbf{q},\mathbf{u}}  \leq \nu_l^{q_l} \sum_{(n_j)_{j\neq l} \in \bN^{r} } \left| \frac{1}{  \sum_{j\neq l} u_j\nu_j n_j} \right|^{q_l}\prod_{j\neq k,l} \frac{1}{n_j^{q_j}} .$$
For any fixed $\mathbf{n} = (n_j)_{j\neq k,l} \in \bN^{r-1}$, the correspondence $$n_k \in \bN \mapsto \sum_{j\neq l} u_j\nu_j n_j \in \bZ$$ is injective, so we can estimate
\begin{align*}
    Y_{k,\bm{\nu},\mathbf{q},\mathbf{u}}  &\leq \nu_l^{q_l}  \sum_{(n_j)_{j\neq k,l} \in \bN^{r-1} } \prod_{j\neq k,l} \frac{1}{n_j^{q_j}} \sum_{n_k\geq 1}\left| \frac{1}{  \sum_{j\neq l} u_j\nu_j n_j} \right|^{q_l}\\
    &\leq \nu_l^{q_l}  \sum_{(n_j)_{j\neq k,l} \in \bN^{r-1} } \prod_{j\neq k,l} \frac{1}{n_j^{q_j}}  \sum_{n_k \in \bZ^{\times}} \frac{1}{|n_k|^{q_l}}\\
    &= 2\nu_l^{q_l} \zeta(q_l) \sum_{(n_j)_{j\neq k,l} \in \bN^{r-1} } \prod_{j\neq k,l} \frac{1}{n_j^{q_j}} \\
    &\leq 2\nu_l^{q_l}\prod_{j\neq k} \zeta(q_j).
\end{align*}
Choosing $\nu_l$ as the minimum of the $\nu_j$ with $j\neq k$, the result follows.\\
\qed
}

As a direct corollary, summing all contributions of $\mathbf{u}\in\{-1,1\}^{r+1}$ yields the following.
\begin{corollary}
    If $\nu_l = \min\{\nu_j : j\neq k \}$, then the following estimate holds,
    $$|\zeta_{k,\bm{\nu},\mathbf{q}}| \leq  2^r(1 + 4\nu_l^{q_l})\prod_{j\neq k}\zeta(q_j),$$
    for all $\mathbf{q}\in \bN_{>1}^{r+1}$.
\end{corollary}
\proof{Indeed, we have $|\zeta_{k,\bm{\nu},\mathbf{q}}| \leq|Y_{k,\bm{\nu},\mathbf{q}}| +|Z_{k,\bm{\nu},\mathbf{q}}|$ with 
$$|Y_{k,\bm{\nu},\mathbf{q}}| \leq \sum_{\mathbf{u}\in \{-1,1 \}^{r+1}} Y_{k,\bm{\nu},\mathbf{q},\mathbf{u}} \leq 2^{r+2} \nu_l^{q_l} \prod_{j\neq k}\zeta(q_j),$$
$$|Z_{k,\bm{\nu},\mathbf{q}}| \leq \frac{1}{2} \sum_{\mathbf{u}\in \{-1,1 \}^{r+1}} Z_{k,\bm{\nu},\mathbf{q},\mathbf{u}} \leq 2^r \prod_{j\neq k} \zeta(q_j),$$
from which the result follows.
\qed}

\begin{remark}\label{emenose}
Some general properties related to the previous sums are as follows. We have
\begin{equation}
    Y_{k,\bm{\nu},\mathbf{q},(1,\ldots ,1)} =  Z_{k,\bm{\nu},\mathbf{q},(1,\ldots ,1)}  = Y_{k,\bm{\nu},\mathbf{q},(-1,\ldots ,-1)} = Z_{k,\bm{\nu},\mathbf{q},(-1,\ldots ,-1)} = 0. 
    \end{equation}
    Moreover, the symmetry by the antipodal map $\mathbf{u} \to -\mathbf{u}$  gives
    \begin{equation}
        Y_{k,\bm{\nu},\mathbf{q},-\mathbf{u}} = Y_{k,\bm{\nu},\mathbf{q},\mathbf{u}},\quad  Z_{k,\bm{\nu},\mathbf{q},-\mathbf{u}} = Z_{k,\bm{\nu},\mathbf{q},\mathbf{u}}.
    \end{equation}
The first identity follows directly from the fact that the equations $\bracket{\mathbf{u} \bm{\nu},\mathbf{n}} = 0$ and $\bracket{\mathbf{u} \bm{\nu}^{k},\mathbf{n}} = 0$ have no nonzero solutions of the same sign, i.e., 
$(\mathbf{u}\bm{\nu})^{\perp}_{\bN}$ and $ (\mathbf{u}\bm{\nu}^{k})^{\perp}_{\bN}$ are empty for $\mathbf{u} = \pm (1,\ldots ,1)$. The second one follows from
$$(\mathbf{u}\bm{\nu})^{\perp}_{\bN} 
 = \{ \mathbf{n} : \bracket{n, \mathbf{u}\bm{\nu}} =0  \} = \{ \mathbf{n} : \bracket{n, -\mathbf{u}\bm{\nu}} =0  \} = (-\mathbf{u}\bm{\nu})^{\perp}_{\bN}.$$
\end{remark}

\subsection{Desingularization Algorithms for Polyhedra}

Coming back to the sums $Y_{k,\bm{\nu},\mathbf{q},\mathbf{u}}$ and $Z_{k,\bm{\nu},\mathbf{q},\mathbf{u}}$, we know they depend on the sets
$(\mathbf{u}\bm{\nu})^{\perp}_{\bN}$ and $(\mathbf{u}\bm{\nu})^{k,\perp}_{\bN}$ which are the positive solutions to linear systems. By Gordan's lemma \cite{alongordan86,coxlittleschenk2011}, these are normal affine semigroups with a finite set of
minimal generators, called a Hilbert basis. As a consequence, $(\mathbf{u}\bm{\nu})^{\perp}_{\bN}$ and 
$(\mathbf{u}\bm{\nu})^{k,\perp}_{\bN}$ are the sets of positive integer points of polyhedral cones $C_{\mathbf{u}}$ and $C^k_{\mathbf{u}}$ respectively with,
$$\dim C_{\mathbf{u}}^k = \dim C_{\mathbf{u}} - 1 = r-1,$$
whose generators lie in  $\bZ^{r+1}$.  For computational purposes, Hilbert bases are not completely satisfactory
since they do not generate the semigroups freely, i.e., representations are not unique, and computing a Hilbert basis becomes increasingly difficult as the dimension $r$ grows.  We can avoid this problem by refining the cones $C_{\mathbf{u}}$ and $C^k_{\mathbf{u}}$ into unimodular cones, i.e., by taking a \emph{desingularization of the cones}. Concretely, there exist finite collections of simplicial cones
$$C_{\mathbf{u},0}, C_{\mathbf{u},1},\ldots ,C_{\mathbf{u},s} \subset  C_{\mathbf{u}},\quad C_{\mathbf{u},0}^k, C_{\mathbf{u},1}^k,\ldots ,C_{\mathbf{u},s_k}^k \subset  C_{\mathbf{u}}^k,$$
whose unions cover the larger cones and which are all of multiplicity one, that is, the primitive generators of each subcone form a $\mathbb{Z}$-basis
of the lattice of integer points in its linear span. These subdivisions define fans $\Sigma_{\mathbf{u}}$ and $\Sigma_{\mathbf{u},k}$, respectively. Denote by $\Sigma^{(d)}_{\mathbf{u}}$ and $\Sigma^{(d)}_{\mathbf{u},k}$ the set of $d$-dimensional cones of each fan. We denote their elements by $F$, and we refer to them as faces of the fans. Denote by $F_{\bN}$ the set of positive integer points of each face in the relative interior of $F$. Then, we have
\begin{equation}
     (\mathbf{u}\bm{\nu})^{\perp}_{\bN} = \bigsqcup_{d=1}^{r} \bigsqcup_{F \in \Sigma^{(d)}_{\mathbf{u}}} F_{\bN},\quad 
    (\mathbf{u}\bm{\nu})^{k,\perp}_{\bN} = \bigsqcup_{d=1}^{r-1} \bigsqcup_{F \in \Sigma^{(d)}_{\mathbf{u},k}} F_{\bN}.
\end{equation}
 Therefore, the previous discussion is summarized as follows.
 \begin{lemma}\label{opensubdiv} Under the previous hypothesis we have
     $$Y_{k,\bm{\nu},\mathbf{q},\mathbf{u}} = \sum_{d=1}^{r } \sum_{F \in \Sigma^{(d)}_{\mathbf{u}}}  \zeta(\mathbf{q}_k , F)  ,\quad Z_{k,\bm{\nu},\mathbf{q},\mathbf{u}} =  \sum_{d=1}^{r-1} \sum_{F \in \Sigma^{(d)}_{\mathbf{u},k}}\zeta(\mathbf{q}_k , F).$$
 \end{lemma}
\qed

This statement is known as an \emph{open subdivision relation} \cite[Lemma 2.7]{conicalzeta2014}. Equivalently, using the notation of \Cref{closedzeta}, we can write
\begin{equation}\label{closedYZ}
Y_{k,\bm{\nu},\mathbf{q},\mathbf{u}}
=
\zeta_{\mathrm{cl}}( \mathbf{q}^k,C_{\mathbf{u}},\Sigma_{\mathbf{u}}),
\qquad
Z_{k,\bm{\nu},\mathbf{q},\mathbf{u}}
=
\zeta_{\mathrm{cl}}( \mathbf{q}^k,C^k_{\mathbf{u}},\Sigma_{\mathbf{u},k}).
\end{equation}

\section{The Case of Dimension $2$}\label{dimension2}  
 In what follows, we will work in the case $r=2$ using the previous results. From \Cref{emenose} we know that, 
$$Y_{k,\bm{\nu},\mathbf{q},\mathbf{u}} = Y_{k,\bm{\nu},\mathbf{q},-\mathbf{u}},\quad Z_{k,\bm{\nu},\mathbf{q},\mathbf{u}} = Z_{k,\bm{\nu},\mathbf{q},-\mathbf{u}} $$
$$Y_{k,\bm{\nu},\mathbf{q},(1,1,1)} = Z_{k,\bm{\nu},\mathbf{q},(1,1,1)}  =  0,$$
for $\mathbf{u}\in \{-1,1\}^3$.
Therefore, our computations are reduced to the vectors
$$\mathbf{u}_0 = (1,-1,-1), \mathbf{u}_1=(1,-1,1),\mathbf{u}_2=(1,1,-1),$$
with 
$$Z_{k,\bm{\nu},\mathbf{q}} =  \sum_{l=0}^2 \sigma_{k,l}Z_{k,\bm{\nu},\mathbf{q},\mathbf{u}_l},\quad Y_{k,\bm{\nu},\mathbf{q}} = 2 \sum_{l=0}^2 \sigma_{k,l}Y_{k,\bm{\nu},\mathbf{q},\mathbf{u}_l},\quad \sigma_{k,l} := \sigma_k(\mathbf{u}_l).$$
Since the terms with $|\mathbf{q}_k|$ odd do not contribute, we only need
to consider the case where $|\mathbf{q}_k|$ is even, for which $\sigma_{k,k} = 1$. As a first step, denote by $C(\mathbf{v}_j,\mathbf{v}_k) = \bR_{\geq 0}\mathbf{v}_j +\bR_{\geq 0}\mathbf{v}_k$ the cone generated by two vectors $\mathbf{v}_j,\mathbf{v}_k \in \bR^3$. Observe that we have the following cone descriptions:
$$(\mathbf{u}_0\bm{\nu})^{\perp}_{\bN} =C(\mathbf{v}_1,\mathbf{v}_2)\cap \bN^{3},$$
$$(\mathbf{u}_1\bm{\nu})^{\perp}_{\bN} = C(\mathbf{v}_2,\mathbf{v}_0)\cap \bN^{3},$$
$$(\mathbf{u}_2\bm{\nu})^{\perp}_{\bN} = C(\mathbf{v}_0,\mathbf{v}_1)\cap \bN^{3},$$
where 
$$\mathbf{v}_0 = \begin{pmatrix}
    0\\
    \nu_2\\
    \nu_1
\end{pmatrix},\quad \mathbf{v}_1 = \begin{pmatrix}
    \nu_2\\
    0\\
    \nu_0
\end{pmatrix},\quad \mathbf{v}_2 = \begin{pmatrix}
     \nu_1\\
    \nu_0\\
   0
\end{pmatrix}. $$
The generators are ordered as in the preceding cone descriptions. In this case, we can compute directly
 $$\bm{\nu}^{k,\perp}_{\bZ} = \bZ \mathbf{u}_l \mathbf{v}_k,\quad l\neq k,$$ 
so, analogously to \Cref{zetaofvectors} we get the following.

\begin{lemma}\label{calculoZr=2} We have
    $$ Z_{k,\bm{\nu},\mathbf{q}} =  \dfrac{\zeta(|\mathbf{q}_k|)}{\mathbf{v_k}^{\mathbf{q}_k}} \sum_{l\neq k}\sigma_{k,l}.$$
\end{lemma}
\qed

For the value $Y_{k,\bm{\nu},\mathbf{q},\mathbf{u}_l}$ we have $(\mathbf{u}_l\bm{\nu})^{\perp}_{\bN} = C_l\cap \bN^3$  
where the planar cone $C_l = C(\mathbf{v}_{l'}, \mathbf{v}_{l''})$ has multiplicity $\nu_l$. By the Hirzebruch--Jung algorithm, this cone admits an explicit
unimodular subdivision  (see \cite[Ch. 10]{coxlittleschenk2011} or \cite[Sec. 2]{torresnova2023}). Thus, for each $l=0,1,2$ we obtain primitive vectors $\mathbf{v}_{0,l},\ldots ,\mathbf{v}_{s_l+1,l} \in \bN^3$ defining an unimodular fan $\Sigma_l$ subdividing $C_l$. The two-dimensional cones of this fan are 
$$C_{l,i} = C(\mathbf{v}_{l,i},\mathbf{v}_{l,i+1}),\quad i = 0,\ldots,s_l.$$
Note that $\Sigma_l^{(1)}$ is simply the set of ray cones defined by the vectors $\mathbf{v}_{l,i}$. Then, we have the following.

\begin{lemma}
     If we denote,
$$\cQ_{k,\bm{\nu},\mathbf{q}} := \sum_{l=0}^2 \sigma_{k,l}\left(\zeta_{cl}^*( \mathbf{q}^k,C_{l},\Sigma_{l})-\frac{1}{\mathbf{v}_{l,0}^{\mathbf{q}_k}}-\frac{1}{\mathbf{v}_{l,s_{l}+1}^{\mathbf{q}_k}} \right) =  \sum_{l=0}^2 \sigma_{k,l}\left(\sum_{i=1}^{s_l} 
\frac{1}{\mathbf{v}_{l,i}^{\mathbf{q}_k}}  + \sum_{i=0}^{s_l}\zeta^*(\mathbf{q}_{k},C_{l,i})\right),$$
    then we have,
    $$\sum_{q_k = 1}\nu_k Y_{k,\bm{\nu},\mathbf{q}} = 2\zeta(|\overline{\mathbf{q}}|)\sum_{q_k = 1} \nu_k \cQ_{k,\bm{\nu},\mathbf{q}}.$$
\end{lemma}
\proof{Assume that the nonzero coordinates of $\mathbf{q}_k$ are $q_{k'},q_{k''}$. Explicitly  we have  $$Y_{k,\bm{\nu},\mathbf{q}, \mathbf{u}_l} = \sum_{\mathbf{n}\in(\mathbf{u}_l\nu)^{\perp}_{\bN}} \frac{1}{n_{k'}^{q_{k'}}n_{k''}^{q_{k''}}}.$$
By the previous discussion, the sets $(\mathbf{u}_l\bm{\nu})^{\perp}_{\bN}$ can be described as a  union
$$(\mathbf{u}_l\bm{\nu})^{\perp}_{\bN} = \bigcup_{i=0}^{s_l}C_{l,i}\cap \bN^{3}.$$
Since each cone 
$C_{l,i}$ has multiplicity one,  from \Cref{opensubdiv} we have
$$ Y_{k,\bm{\nu},\mathbf{q}, \mathbf{u}_l} =
\sum_{i=0}^{s_l}\zeta(\mathbf{q}_k,C_{l,i}) +   \sum_{i=1}^{s_l} 
\zeta(\mathbf{q}_k, \mathbf{v}_{l,i}  ).$$
From \Cref{zetaofvectors} and the fact that $\zeta(\mathbf{q}_k,C_{l,i}) = \zeta(|\mathbf{q}_k|)\zeta^*(\mathbf{q}_k,C_{l,i})$ we obtain
$$ Y_{k,\bm{\nu},\mathbf{q}, \mathbf{u}_l} = 
\zeta(|\mathbf{q}_k|) \left(\sum_{i=0}^{s_l}\zeta^*(\mathbf{q}_{k},C_{l,i})  +   \sum_{i=1}^{s_l} 
\dfrac{1}{\mathbf{v}_{l,i}^{\mathbf{q}_k}} \right).$$
Therefore, we put everything together with \Cref{convergence} to ensure convergence, replacing $|\mathbf{q}_k| = |\overline{\mathbf{q}}|$ for all $k$ such that $q_k=1$.\\
\qed
}

\begin{corollary}
For all $\mathbf{q}\in \bN^{3}$, the reciprocity formula for Bernoulli polynomials is
$$\cR_{b_{\mathbf{q}}}(\bm{\nu}) = -\frac{ B_{|\overline{\mathbf{q}}|}}{|\overline{\mathbf{q}}|!} \prod_{q_j>1} q_j! \sum_{q_k=1} \nu_k \left(\frac{1}{\mathbf{v}_k^{\mathbf{q}_k}}\sum_{k\neq l} \frac{\sigma_{k,l}}{2} +  \cQ_{k,\bm{\nu},\mathbf{q}}  \right) -  \frac{1 -  (-1)^{\mathbf{1}_{\mathbf{q}}}}{2^{\mathbf{1}_{\mathbf{q}}}}  \prod_{q_k > 1} B_{q_k}   .$$
\end{corollary}
\proof{From \Cref{reciforbern} we know that
$$\cR_{b_{\mathbf{q}}}(\bm{\nu}) = I_{\overline{\mathbf{q}}}\sum_{q_k=1} \nu_k \zeta_{k,\bm{\nu},\mathbf{q}} -  \frac{1 -  (-1)^{\mathbf{1}_{\mathbf{q}}}}{2^{\mathbf{1}_{\mathbf{q}}}}  \prod_{q_k > 1} B_{q_k} ,\quad I_{\overline{\mathbf{q}}} :=  \frac{(-1)^{|\overline{\mathbf{q}}|/2}}{(2\pi)^{|\overline{\mathbf{q}}|}}\prod_{q_j>1} q_j!,$$
for $|\overline{\mathbf{q}}|$ even. The previous two lemmas give 
$$\sum_{q_k=1}\nu_k \zeta_{k,\bm{\nu},\mathbf{q}} =  2\zeta(|\overline{\mathbf{q}}|) \sum_{q_k=1}\nu_k \left(\frac{1}{\mathbf{v}_k^{\mathbf{q}_k}}\sum_{k\neq l} \frac{\sigma_{k,l}}{2} +  \cQ_{k,\bm{\nu},\mathbf{q}}  \right) .$$
Together with the identity
$$\zeta(|\overline{\mathbf{q}}|) = (-1)^{\frac{|\overline{\mathbf{q}}|}{2}+1}\frac{B_{|\overline{\mathbf{q}}|}(2\pi)^{|\overline{\mathbf{q}}|}}{2|\overline{\mathbf{q}}|!},$$
we get the desired result.\\
\qed
}

\begin{remark}
    In the case where $\mathbf{q}$ satisfies $q_k = 1$ and $q_{k'} = q_{k''} = q$ for $q\geq 1$ we have $\sigma_{k,l} = (-1)^q$ for all $k\neq l$. Then, the previous identity reduces to
   $$\cR_{b_{\mathbf{q}}}(\bm{\nu}) = S_{b_\mathbf{q}}(\bm{\nu}|\nu_k) =  -\frac{B_{2q}q!^2 \nu_k}{(2q)!} \left(\frac{(-1)^q}{\nu_{k'}^q\nu_{k''}^q} +  \cQ_{k,\bm{\nu},\mathbf{q}}  \right)- B_q^2,$$
   for $q>1$. From \Cref{vanishremark}, when $q = 1$, we know that $\sum_{k=0}^2\nu_k \cQ_{k,\bm{\nu},\mathbf{q}} = 0$, so our expression simplifies to
   $$\cR_{b_{\mathbf{q}}}(\bm{\nu}) = \frac{\nu_0^2 + \nu_1^2+\nu_2^2}{12\nu_0 \nu_1\nu_2} - \frac{1}{4}$$
   which is just Dedekind-Rademacher's reciprocity. In the case $q>1$, applying the previous identity to the three cyclic choices
$
(1,q,q),(q,1,q),$ and $(q,q,1),
$
and denoting $\cQ_{\bm{\nu},\mathbf{q}} = \sum_{k=0}^2\nu_k \cQ_{k,\bm{\nu},\mathbf{q}}$,  we obtain an analogue of
Dedekind--Rademacher's reciprocity: 
$$ \sum_{k=0}^2 S_{b_\mathbf{q}}(\bm{\nu}|\nu_k) = -\frac{B_{2q}q!^2}{(2q)!} \left(\frac{\nu_0^{q+1}+\nu_1^{q+1} +\nu_2^{q+1}}{\nu_0^q\nu_1^q \nu_2^q}(-1)^q +   \cQ_{\bm{\nu},\mathbf{q}}  \right)-  3B_q^2.$$ 
This shows how sums involving continuous Bernoulli factors can be studied
by inserting one \(b_1\)-factor and using the corresponding single-jump
reciprocity.
\end{remark}

\begin{example}
As a concrete instance of the preceding remark, let
$$
    \bm{\nu}=(2,a,b),
    \qquad 
    \mathbf{q}=(1,q,q),
$$
where $a,b\geq 3$ are odd, coprime, and $a\neq b$. Assume also that
$q\geq 2$. Since the only component of $\mathbf q$ equal to $1$ is
$q_0$, the reciprocity symbol reduces to the single Dedekind sum
$$
    R_{b_{\mathbf{q}}}(\bm{\nu})
    =
    S_{b_{\mathbf{q}}}(\bm{\nu} | 2).
$$
Explicitly, since $a$ and $b$ are odd, we have
$$
    S_{b_{\mathbf{q}}}((2,a,b)|2)
    =
    b_q\left(\frac{a}{2}\right)
    b_q\left(\frac{b}{2}\right) =
    b_q\left(\frac12\right)^2.
$$
Using the identity $
    b_q\left(\frac12\right)
    =
    (2^{1-q}-1)B_q,$
we get
$$
    S_{b_{\mathbf{q}}}((2,a,b) |2)
    =
    (2^{2-2q} - 2^{2-q}+1)B_q^2.
$$
On the conical side, the vectors associated with
\(\bm{\nu}=(2,a,b)\) are
$$
    \mathbf{v}_0=(0,b,a),\quad
    \mathbf{v}_1=(b,0,2),\quad
    \mathbf{v}_2=(a,2,0),
$$
defining the three relevant cones
$$
    C(\mathbf{v}_1,\mathbf{v}_2),\quad C(\mathbf{v}_2,\mathbf{v}_0),\quad C(\mathbf{v}_0,\mathbf{v}_1),
$$
with multiplicities $2,a,b$ respectively. By the previous remark, applied to $k=0$, we have
$$
    S_{b_{\mathbf{q}}}((2,a,b) | 2)
    =
    -\frac{2B_{2q}(q!)^2}{(2q)!}
    \left(
        \frac{(-1)^q}{a^q b^q}
        +
        Q_{0,(2,a,b),\mathbf{q}}
    \right)
    -
    B_q^2.
$$
Hence,
$$
    Q_{0,(2,a,b),\mathbf{q}}
    =
    -\frac{(2q)!B_q^2}{B_{2q}(q!)^2}
    \left(
       1 - 2^{1-q} + 2^{1-2q}
    \right)
    -
    \frac{(-1)^q}{a^q b^q}.
$$
Let us now isolate the contribution of the cone of multiplicity $2$.
For $k=0$, we have $
    \mathbf{q}^0=(0,q,q)$,
and
$$
    Q_{0,(2,a,b),\mathbf{q}}
    =
    \mathcal A_0^{(q)}
    +
    (-1)^q\mathcal A_1^{(q)}
    +
    (-1)^q\mathcal A_2^{(q)},
$$
$$\mathcal A_l^{(q)} = \zeta_{cl}^*( \mathbf{q}^0,C_{l},\Sigma_{l}) - \frac{1}{\mathbf{v}_{l,0}^{\mathbf{q}^0}}- \frac{1}{\mathbf{v}_{l,s_l+1}^{\mathbf{q}^0}} = \sum_{i=1}^{s_l} 
\frac{1}{\mathbf{v}_{l,i}^{\mathbf{q}^0}}  + \sum_{i=0}^{s_l}\zeta^*(\mathbf{q}^0,C_{l,i}),$$
coming from the
Hirzebruch--Jung subdivision of the $l$-th cone. The first term
$\mathcal{A}_0^{(q)}$ corresponds to
$$
    C_0 = C(\mathbf{v}_1,\mathbf{v}_2)=C((b,0,2),(a,2,0)).
$$
Since this cone has multiplicity $2$, its subdivision has a single
interior ray, generated by
$$
    \mathbf{w}=\frac{\mathbf{v}_1+\mathbf{v}_2}{2}
    =
    \left(\frac{a+b}{2},1,1\right).
$$
Thus $C(\mathbf{v}_1,\mathbf{v}_2)
    = C_{0,0} \cup C_{0,1}$, where $C_{0,0} = C(\mathbf{v}_1,\mathbf{w})$ and $C_{0,1} =  C(\mathbf{w},\mathbf{v}_2),$
 both unimodular cones. Therefore
$$
    \mathcal A_0^{(q)} = \zeta_{cl}^*( \mathbf{q}^0,C_{0},\Sigma_{0})
    =
    1
    +
    \zeta^*(\mathbf{q}^0,C_{0,0})
    +
    \zeta^*(\mathbf{q}^0,C_{0,1}).
$$
By symmetry we have 
$\zeta^*(\mathbf{q}^0,C_{0,0})
    =
    \zeta^*(\mathbf{q}^0,C_{0,1})$, and explicitly
    $$\zeta^*(\mathbf{q}^0,C_{0,0}) =\frac{1}{\zeta(2q)}\zeta(\mathbf{q}^0,C_{0,0}) =
    \frac{1}{\zeta(2q)}
    \sum_{m,n\geq 1}
    \frac{1}{n^q(2m+n)^q}.
$$
Now we compute,
\begin{align*}
    \sum_{m,n\geq 1}
    \frac{1}{n^q(2m+n)^q} &= \sum_{\substack{1\leq n < N\\ n,N \text{ even}} } \frac{1}{n^qN^q}+\sum_{\substack{1\leq n < N\\ n,N \text{ odd}} } \frac{1}{n^qN^q}\\
    & = \sum_{\substack{1\leq n < N} } \frac{2^{-2q}}{n^qN^q}+\sum_{\substack{1\leq n < N }} \frac{1}{(2n-1)^q(2N-1)^q}\\
    &= 2^{-2q}\left(\frac{\zeta(q)^2 - \zeta(2q)}{2}\right) + \frac{1}{2}\left[\left( \sum_{n\geq 1} \frac{1}{(2n-1)^q}\right)^2 - \sum_{n\geq 1} \frac{1}{(2n-1)^{2q}}\right]\\
    &= 2^{-2q}\left(\frac{\zeta(q)^2 - \zeta(2q)}{2}\right) + \frac{1}{2}\left[(1-2^{-q})^2\zeta(q)^2 - (1-2^{-2q})\zeta(2q)\right]\\
    &= \frac{1}{2}\left[ (1-2^{1-q} + 2^{1-2q})\zeta(q)^2 - \zeta(2q)\right].
\end{align*}
So, we get
$$
    \zeta^*(\mathbf{q}^0,C_{0,0}))
    =
    \frac{1}{2}
    \left[
        \left(
            1-2^{1-q}+2^{1-2q}
        \right)
        \frac{\zeta(q)^2}{\zeta(2q)}
        -1
    \right],
$$
and consequently,
$$
    \mathcal A_0^{(q)}
    =
    \left(
        1-2^{1-q}+2^{1-2q}
    \right)
    \frac{\zeta(q)^2}{\zeta(2q)}.
$$
Thus the reciprocity formula determines the remaining contribution
coming from the cones of multiplicities $a$ and $b$. Namely,
$$   
        \mathcal A_1^{(q)}+\mathcal A_2^{(q)}
    =
   (-1)^q \left(
        2^{1-q} - 2^{1-2q} - 1
    \right)\left( \frac{B_q^2 (2q)!}{B_{2q}(q!)^2} + \frac{\zeta(q)^2}{\zeta(2q)} \right)  
    -
    \frac{1}{a^q b^q},
$$
or equivalently
$$   
      \zeta_{cl}^*( \mathbf{q}^0,C_{1},\Sigma_{1})+\zeta_{cl}^*( \mathbf{q}^0,C_{2},\Sigma_{2})
    =
   (-1)^q \left(
        2^{1-q} - 2^{1-2q} - 1
    \right)\left( \frac{B_q^2 (2q)!}{B_{2q}(q!)^2} + \frac{\zeta(q)^2}{\zeta(2q)} \right)  
    +
    \frac{1}{a^q b^q}.
$$
This identity shows that, although the individual reduced conical zeta
values arising from the cones of multiplicities $a$ and $b$ are not
evaluated separately, their total contribution is explicitly determined
by the reciprocity formula.
\end{example}

\bibliographystyle{plainnat}
\bibliography{ref.bib}
\vspace{5mm}

\end{document}